\documentclass{article}
\usepackage[top=30truemm,bottom=30truemm,left=30truemm,right=30truemm]{geometry}
\usepackage{amssymb,amsmath,amsthm,amscd,ascmac}
\usepackage{bm}
\usepackage{tikz}
\usetikzlibrary{intersections,calc,arrows.meta}
\newtheorem{theorem}{Theorem}

\newtheorem*{example}{Example}
\theoremstyle{definition}

\begin{document}
\title{Repeated Eigenvalues Imply Nodes? \\
A Problem of Planar Differential Equations}
\markright{Do Repeated Eigenvalues Imply Nodes?}
\author{Kenzi Odani\\
Department of mathematics, Aichi University of Education}
\maketitle
\begin{abstract}
Poincar\'e gave a criterion which determines the shape 
of equilibrium for planar differential equations. 
In his statement, he excluded the case of repeated eigenvalues. 
In fact, in such a case, we can give a $C^1$ counter-example to 
his assertion. 
In this note, we show that if we strengthen the condition to 
$C^{1,\alpha}$ ($0<\alpha<1$), his assertion becomes true 
even in case of repeated eigenvalues.
\end{abstract}
\section{Introduction.}
Consider a planar differential equation below: 
\begin{align}
\dfrac{d\bm{x}}{dt}=\bm{f}(\bm{x}),\quad\bm{x}\in\mathbb{R}^2, 
\label{eqn:1}
\end{align}
where $\bm{f}$ is a 2-dimensional function. 
We call a point $\bm{p}$ an {\it equilibrium} if $\bm{f}(\bm{p})=\bm{0}$. 
We call an equilibrium {\it repelling (attracting)} if all 
nearby trajectories are outgoing (incoming). 
We call a repelling or an attracting equilibrium $\bm{p}$ 
a {\it node (focus)} if the direction 
$\bm{d}(t)=\bm{f}\big(\bm{x}(t)\big)/|\bm{f}\big(\bm{x}(t)\big)|$
converges to a constant vector (rotates infinite times) 
as $\bm{x}(t)\to\bm{p}$. 
We call an equilibrium a {\it saddle} if it has two outgoing and 
two incoming trajectories, which appear alternately. 
We give their phase portraits below. \par
\medskip
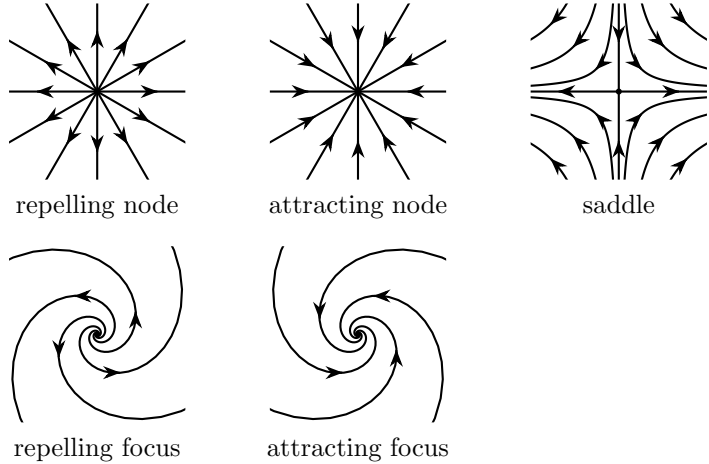
\begin{center}
\begin{tabular}{c@{\hspace{3em}}c@{\hspace{3em}}c}
\begin{tikzpicture}[samples=300,scale=1.00]
\begin{scope}\clip(-1.2,-1.2) rectangle (1.2,1.2);
\draw[thick](1.697,0)--(-1.697,0);
\draw[thick](1.47,0.849)--(-1.47,-0.849);
\draw[thick](0.849,1.47)--(-0.849,-1.47);
\draw[thick](0,1.697)--(0,-1.697);
\draw[thick](-0.849,1.47)--(0.849,-1.47);
\draw[thick](-1.47,0.849)--(1.47,-0.849);
\draw[-Stealth,thick](0.8,0)--(0.81,0);
\draw[-Stealth,thick](0.693,0.4)--(0.701,0.405);
\draw[-Stealth,thick](0.4,0.693)--(0.405,0.701);
\draw[-Stealth,thick](0,0.8)--(0,0.81);
\draw[-Stealth,thick](-0.4,0.693)--(-0.405,0.701);
\draw[-Stealth,thick](-0.693,0.4)--(-0.701,0.405);
\draw[-Stealth,thick](-0.8,0)--(-0.81,0);
\draw[-Stealth,thick](-0.693,-0.4)--(-0.701,-0.405);
\draw[-Stealth,thick](-0.4,-0.693)--(-0.405,-0.701);
\draw[-Stealth,thick](0,-0.8)--(0,-0.81);
\draw[-Stealth,thick](0.4,-0.693)--(0.405,-0.701);
\draw[-Stealth,thick](0.693,-0.4)--(0.701,-0.405);
\fill[black](0,0)circle(0.04);
\draw[ultra thick,white](1.2,1.2)--(-1.2,1.2)--(-1.2,-1.2)--(1.2,-1.2)--cycle;
\end{scope}
\end{tikzpicture} & 
\begin{tikzpicture}[samples=300,scale=1.00]
\begin{scope}\clip(-1.2,-1.2) rectangle (1.2,1.2);
\draw[thick](1.697,0)--(-1.697,0);
\draw[thick](1.47,0.849)--(-1.47,-0.849);
\draw[thick](0.849,1.47)--(-0.849,-1.47);
\draw[thick](0,1.697)--(0,-1.697);
\draw[thick](-0.849,1.47)--(0.849,-1.47);
\draw[thick](-1.47,0.849)--(1.47,-0.849);
\draw[-Stealth,thick](0.66,0)--(0.65,0);
\draw[-Stealth,thick](0.572,0.33)--(0.563,0.325);
\draw[-Stealth,thick](0.33,0.572)--(0.325,0.563);
\draw[-Stealth,thick](0,0.66)--(0,0.65);
\draw[-Stealth,thick](-0.33,0.572)--(-0.325,0.563);
\draw[-Stealth,thick](-0.572,0.33)--(-0.563,0.325);
\draw[-Stealth,thick](-0.66,0)--(-0.65,0);
\draw[-Stealth,thick](-0.572,-0.33)--(-0.563,-0.325);
\draw[-Stealth,thick](-0.33,-0.572)--(-0.325,-0.563);
\draw[-Stealth,thick](0,-0.66)--(0,-0.65);
\draw[-Stealth,thick](0.33,-0.572)--(0.325,-0.563);
\draw[-Stealth,thick](0.572,-0.33)--(0.563,-0.325);
\fill[black](0,0)circle(0.04);
\draw[ultra thick,white](1.2,1.2)--(-1.2,1.2)--(-1.2,-1.2)--(1.2,-1.2)--cycle;
\end{scope}
\end{tikzpicture} & 
\begin{tikzpicture}[samples=300,scale=1.00]
\begin{scope}\clip(-1.2,-1.2) rectangle (1.2,1.2);
\draw[thick](1.2,0)--(-1.2,0);
\draw[thick](0,1.2)--(0,-1.2,0);
\draw[thick,variable=\t,domain=0.01:1.5,samples=300]
  plot ({\t},{0.09/\t});
\draw[thick,variable=\t,domain=0.01:1.5,samples=300]
  plot ({\t},{0.36/\t});
\draw[thick,variable=\t,domain=0.01:1.5,samples=300]
  plot ({\t},{0.81/\t});
\draw[thick,variable=\t,domain=0.01:1.5,samples=300]
  plot ({-\t},{-0.09/\t});
\draw[thick,variable=\t,domain=0.01:1.5,samples=300]
  plot ({-\t},{-0.36/\t});
\draw[thick,variable=\t,domain=0.01:1.5,samples=300]
  plot ({-\t},{-0.81/\t});
\draw[thick,variable=\t,domain=0.01:1.5,samples=300]
  plot ({-\t},{0.09/\t});
\draw[thick,variable=\t,domain=0.01:1.5,samples=300]
  plot ({-\t},{0.36/\t});
\draw[thick,variable=\t,domain=0.01:1.5,samples=300]
  plot ({-\t},{0.81/\t});
\draw[thick,variable=\t,domain=0.01:1.5,samples=300]
  plot ({\t},{-0.09/\t});
\draw[thick,variable=\t,domain=0.01:1.5,samples=300]
  plot ({\t},{-0.36/\t});
\draw[thick,variable=\t,domain=0.01:1.5,samples=300]
  plot ({\t},{-0.81/\t});
\draw[-Stealth,thick](0.8,0)--(0.81,0);
\draw[-Stealth,thick](-0.8,0)--(-0.81,0);
\draw[-Stealth,thick](0,0.66)--(0,0.65);
\draw[-Stealth,thick](0,-0.66)--(0,-0.65);
\draw[-Stealth,thick](0.3,0.3)--(0.31,0.29);
\draw[-Stealth,thick](0.6,0.6)--(0.61,0.59);
\draw[-Stealth,thick](0.9,0.9)--(0.91,0.89);
\draw[-Stealth,thick](-0.3,0.3)--(-0.31,0.29);
\draw[-Stealth,thick](-0.6,0.6)--(-0.61,0.59);
\draw[-Stealth,thick](-0.9,0.9)--(-0.91,0.89);
\draw[-Stealth,thick](-0.3,-0.3)--(-0.31,-0.29);
\draw[-Stealth,thick](-0.6,-0.6)--(-0.61,-0.59);
\draw[-Stealth,thick](-0.9,-0.9)--(-0.91,-0.89);
\draw[-Stealth,thick](0.3,-0.3)--(0.31,-0.29);
\draw[-Stealth,thick](0.6,-0.6)--(0.61,-0.59);
\draw[-Stealth,thick](0.9,-0.9)--(0.91,-0.89);
\fill[black](0,0)circle(0.04);
\draw[ultra thick,white](1.2,1.2)--(-1.2,1.2)--(-1.2,-1.2)--(1.2,-1.2)--cycle;
\end{scope}
\end{tikzpicture} \\
repelling node & attracting node & saddle \\[0.8em]
\begin{tikzpicture}[samples=300,scale=1.00]
\begin{scope}\clip(-1.2,-1.2) rectangle (1.2,1.2);
\draw[thick,variable=\t,domain=-1:50,samples=300]
  plot ({exp(-\t*0.5)*cos(\t r)},{-exp(-\t*0.5)*sin(\t r)});
\draw[thick,variable=\t,domain=-1:50,samples=300]
  plot ({-exp(-\t*0.5)*cos(\t r)},{exp(-\t*0.5)*sin(\t r)});
\draw[thick,variable=\t,domain=-1:50,samples=300]
  plot ({exp(-\t*0.5)*sin(\t r)},{exp(-\t*0.5)*cos(\t r)});
\draw[thick,variable=\t,domain=-1:50,samples=300]
  plot ({-exp(-\t*0.5)*sin(\t r)},{-exp(-\t*0.5)*cos(\t r)});
\draw[-Stealth,thick](-0.30,0.51)--(-0.31,0.51);
\draw[-Stealth,thick](-0.51,-0.30)--(-0.51,-0.31);
\draw[-Stealth,thick](0.30,-0.51)--(0.31,-0.51);
\draw[-Stealth,thick](0.51,0.30)--(0.51,0.31);
\fill[black](0,0)circle(0.04);
\draw[ultra thick,white](1.2,1.2)--(-1.2,1.2)--(-1.2,-1.2)--(1.2,-1.2)--cycle;
\end{scope}
\end{tikzpicture} & 
\begin{tikzpicture}[samples=300,scale=1.00]
\begin{scope}\clip(-1.2,-1.2) rectangle (1.2,1.2);
\draw[thick,variable=\t,domain=-1:50,samples=300]
  plot ({exp(-\t*0.5)*cos(\t r)},{exp(-\t*0.5)*sin(\t r)});
\draw[thick,variable=\t,domain=-1:50,samples=300]
  plot ({-exp(-\t*0.5)*cos(\t r)},{-exp(-\t*0.5)*sin(\t r)});
\draw[thick,variable=\t,domain=-1:50,samples=300]
  plot ({-exp(-\t*0.5)*sin(\t r)},{exp(-\t*0.5)*cos(\t r)});
\draw[thick,variable=\t,domain=-1:50,samples=300]
  plot ({exp(-\t*0.5)*sin(\t r)},{-exp(-\t*0.5)*cos(\t r)});
\draw[-Stealth,thick](0.21,0.51)--(0.20,0.51);
\draw[-Stealth,thick](-0.51,0.21)--(-0.51,0.20);
\draw[-Stealth,thick](-0.21,-0.51)--(-0.20,-0.51);
\draw[-Stealth,thick](0.51,-0.21)--(0.51,-0.20);
\fill[black](0,0)circle(0.04);
\draw[ultra thick,white](1.2,1.2)--(-1.2,1.2)--(-1.2,-1.2)--(1.2,-1.2)--cycle;
\end{scope}
\end{tikzpicture} \\
repelling focus & attracting focus 
\end{tabular}
\end{center}\par
\medskip
The notions such as node, saddle and focus were introduced by 
H.~Poincar\'e \cite{P}. 
He gave a criterion which determines the shape of equilibrium. 
Let $\bm{p}$ be an equilibrium of Eq.(\ref{eqn:1}), 
and $\lambda_1,\lambda_2$ the eigenvalues of Jacobian matrix. 
He stated that when $\bm{f}$ is a polynomial vector and 
$\lambda_1\not=\lambda_2$, the following assertions hold: 
\begin{enumerate}
\item[(i)] If $\lambda_1,\lambda_2$ are real with the same sign, 
then $\bm{p}$ is a node.
\item[(ii)] If $\lambda_1,\lambda_2$ are real with opposite signs, 
then $\bm{p}$ is a saddle.
\item[(iii)] If $\lambda_1,\lambda_2$ are imaginary with 
nonzero real part, then $\bm{p}$ is a focus.
\end{enumerate}\par
\noindent
In case of linear systems, we can prove his statement 
even when $\lambda_1=\lambda_2$. The proof can be found in many 
textbooks. See \cite[\S20.4]{A1}, for example. \par
\medskip
%
In case of nonlinear systems, he did not provide any complete proofs. 
However, we can prove his statement for nonlinear 
$C^1$ systems. \par
\begin{theorem}
Suppose that Eq.(\ref{eqn:1}) is of class $C^1$. 
If $\lambda_1\not=\lambda_2$, then the assertions 
(i), (ii) and (iii) hold.
\end{theorem}
The condition $\lambda_1\not=\lambda_2$ is critical. 
In fact, if it is dropped, the following counter-example appears. \par
\begin{example}
Consider a $C^1$ nonlinear system below: 
\begin{align}
\dfrac{d}{dt}\begin{bmatrix}\ x\ \\ \ y\ \end{bmatrix}
=\begin{bmatrix}\ -1 \!&\! 0\ \\ \ \varepsilon \!&\! -1\ \end{bmatrix}
\begin{bmatrix}\ x\ \\ \ y\ \end{bmatrix}
-\dfrac{2}{\ln(x^2+y^2)}
\begin{bmatrix}\,-y\,\\ \,x\,\end{bmatrix}, 
\label{eqn:2}
\end{align}
where $0\leqq\varepsilon<2$. Then the equilibrium $\bm{0}$ is a focus. 
\end{example}
The above example shows that if we drop the condition 
$\lambda_1\not=\lambda_2$, Theorem~1 becomes false. 
However, such an example is very fragile. 
If we slightly strengthen the $C^1$ condition, Theorem~1 
becomes true even when $\lambda_1=\lambda_2$. \par
\begin{theorem}
Suppose that Eq.(\ref{eqn:1}) is of class $C^{1,\alpha}$ ($\alpha>0$). 
Then the assertions (i), (ii) and (iii) hold.
\end{theorem}
Here, we say that the system (\ref{eqn:1}) is of class $C^{1,\alpha}$ 
if $\bm{f}$ are of class $C^1$ and their first partial derivatives 
are $\alpha$-H\"older continuous. 
The $C^{1,\alpha}$ condition ($0<\alpha<1$) is stronger than 
$C^1$ and weaker than $C^2$. \par
\section{Proofs.}
Regarding the proof of Theorem~1, the assertion (ii) is guaranteed 
by Grobman-Hart\-man Theorem, which is found in \cite[\S5]{R}, 
for example. 
Unfortunately, the theorem can not distiguish (i) and (iii) because 
nodes and foci are topologically equivalent. See \cite[\S22.6]{A1}.
Instead, we can prove (i) and (iii) by diagonalizing Jacobian matrix 
and by transforming it to a form of polar coordinates. 
We can carry out the proof without any difficulties. 
So we leave the proof to the reader. \par
\medskip
In this section, we shall prove Example and Theorem 2. 
\begin{proof}[Proof of Example]
We transform Eq.(\ref{eqn:2}) to a form of polar coordinates:
\begin{align}
r'(t)=-r+\varepsilon r\cos\theta\sin\theta,\quad 
\theta'(t)=\varepsilon\cos^2\theta-\dfrac{1}{\ln r}, 
\label{eqn:3}
\end{align}
By evaluating the first of Eq.(\ref{eqn:3}), we have that
\begin{align*}
-k_1r\leqq r'(t)\leqq -k_2r,\quad 
k_1:=1+\tfrac{1}{2}\varepsilon, \ 
k_2:=1-\tfrac{1}{2}\varepsilon.
\end{align*}
Take $r_0=r(0)$ so small that $r_0\in[0,1)$. 
By integrating the above, we have that
\begin{align}
r_0e^{-k_1t}\leqq r(t)\leqq r_0e^{-k_2t}<1.
\label{eqn:4}
\end{align}
Therefore, $r(t)\to 0$ as $t\to\infty$. 
Hence $\bm{0}$ is attracting. \par
\medskip
By taking logarithm of Eq.(\ref{eqn:4}), we have that 
\begin{align}
-(k_1t-\ln r_0)\leqq\ln r(t)<0.
\label{eqn:4.5}
\end{align}
By applying Eq.(\ref{eqn:4.5}) to the second of Eq.(\ref{eqn:3}), 
we have that
\begin{align*}
\theta'(t)\geqq-\dfrac{1}{\ln r(t)}\geqq\dfrac{1}{k_1t-\ln r_0}>0, 
\end{align*}
By integrating it, we have that
\begin{align*}
\theta(t)-\theta(0)\geqq
\dfrac{1}{k_1}\big\{\ln(k_1t-\ln r_0)-\ln(-\ln r_0)\big\}.
\end{align*}
Therefore, $\theta(t)\to\infty$ as $t\to\infty$. 
Hence $\bm{0}$ is a focus.
\end{proof}
\begin{proof}[Proof of Theorem 2.]
By translating $\bm{p}$ to $\bm{0}$, we make Eq.(\ref{eqn:1}) 
the following form:
\begin{align}
\dfrac{d}{dt}\begin{bmatrix}\ x\ \\ \ y\ \end{bmatrix}
=\begin{bmatrix}\ a_{11} \!&\! a_{12}\ \\
\ a_{21} \!&\! a_{22}\ \end{bmatrix}
\begin{bmatrix}\ x\ \\ \ y\ \end{bmatrix}
+\begin{bmatrix}\ u(x,y)\ \\ \ v(x,y)\ \end{bmatrix}, 
\label{eqn:5}
\end{align}
where $u$, $v$ are nonlinear terms.  
Remark that the square matrix of Eq.(\ref{eqn:5}) is Jacobian matrix 
and denote it by $J$. \par
\medskip
Consider only the case $\lambda_1=\lambda_2$. This is because other 
cases are already established in Theorem 1. 
Put $\lambda=\lambda_1$, and assume that $\lambda<0$. 
(It does not hurt the generality.) 
We devide the proof into two cases. \par
\medskip
\textbf{1.} The case when $J$ is diagonalizable. 
In such a case, there is a matrix $P$ such that 
$P^{-1}JP=\Big[\begin{matrix}\ \lambda \!&\! 0\ \\[-0.1em]
\ 0 \!&\! \lambda \ \end{matrix}\Big]$. 
So we can transform Eq.(\ref{eqn:5}) into the following form:
\begin{align*}
\dfrac{d}{dt}\begin{bmatrix}\ x\ \\ \ y\ \end{bmatrix}
=\begin{bmatrix}\ \lambda \!&\! 0\ \\ \ 0 \!&\! \lambda\ \end{bmatrix}
\begin{bmatrix}\ x\ \\ \ y \ \end{bmatrix}
+\begin{bmatrix}\ p(x,y)\ \\ \ q(x,y)\ \end{bmatrix}.
\end{align*}
Also, we transform it to a form of polar coordinates:
\begin{align}
r'(t)=\lambda r+r\varphi(r,\theta),\quad 
\theta'(t)=\psi(r,\theta),
\label{eqn:6}
\end{align}
Since both $\varphi,\psi$ are $\alpha$-H\"older continuous, 
there are reals $c>0,m>0$ such that
\begin{align}
\big|\,\varphi(r,\theta)\,\big|<mr^{\alpha},\quad
\big|\,\psi(r,\theta)\,\big|<mr^{\alpha}\quad
\textrm{on }[0,c]\times\mathbb{R}.
\label{eqn:7}
\end{align}
We take a real $d$ so small that $md^{\alpha}<\tfrac{1}{2}|\lambda|$ 
and $d\leqq c$. 
If we take $r(t)\in[0,d]$, then we evaluate the first of 
Eq.(\ref{eqn:6}) as follows:
\begin{align*}
r'(t)<\big(\!-\!|\lambda|+\tfrac{1}{2}|\lambda|\,\big)r=-kr,\quad 
k:=\tfrac{1}{2}|\lambda|.
\end{align*}
If we take $r_0=r(0)$ so small that $r_0\in[0,d]$, 
then we have that
\begin{align}
r(t)<r_0e^{-kt}.
\label{eqn:8}
\end{align} 
Therefore, $r(t)\to 0$ as $t\to\infty$. 
Hence $\bm{p}$ is attracting. \par
\medskip
By applying Eqs.(\ref{eqn:7}), (\ref{eqn:8}) to the second 
of Eq.(\ref{eqn:6}), we have that
\begin{align}
\big|\,\theta'(t)\,\big|<m\big(r_0e^{-kt}\big)^{\alpha}
=mr_0^{\,\alpha}e^{-lt},\quad 
l:=\alpha k. 
\label{eqn:9}
\end{align}
By integrating Eq.(\ref{eqn:9}) from $0$ to $t$ ($0\leqq t$), 
we have that
\begin{align*}
\big|\,\theta(t)-\theta(0)\,\big|
\leqq-\dfrac{m}{l}r_0^{\,\alpha}\Big(e^{-lt}-1\Big)
<\dfrac{m}{l}r_0^{\,\alpha}. 
\end{align*}
Since $\theta(t)$ is bounded, 
there is a sequence $\{t_n\}$ diverging to $\infty$ such that 
$\theta(t_n)$ converges to a limit $\theta^*$. 
By integrating Eq.(\ref{eqn:9}) from $t$ to $t_n$ ($t\leqq t_n$), 
we have that
\begin{align*}
\big|\,\theta(t_n)-\theta(t)\,\big|
\leqq-\dfrac{m}{l}r_0^{\,\alpha}\Big(e^{-lt_n}-e^{-lt}\Big)
<\dfrac{m}{l}r_0^{\,\alpha}e^{-lt}. 
\end{align*}
By taking $n\to\infty$, we have that
\begin{align*}
\big|\,\theta^*-\theta(t)\,\big|\leqq\dfrac{m}{l}r_0^{\,\alpha}e^{-lt}. 
\end{align*}
Therefore, $\theta(t)$ converges to $\theta^*$. 
Hence $\bm{p}$ is a node. \par
\medskip
\textbf{2.} The case when $J$ is not diagonalizable. 
In such a caes, there is a matrix $P$ such that 
$P^{-1}JP=\Big[\begin{matrix}\ \lambda \!&\! 0\ \\[-0.1em]
\ \varepsilon \!&\! \lambda\ \end{matrix}\Big]$ (Jordan normal form). 
Put $\varepsilon=-\lambda$.
So we can transform Eq.(\ref{eqn:5}) into the following form:
\begin{align*}
\dfrac{d}{dt}\begin{bmatrix}\ x\ \\ \ y\ \end{bmatrix}
=\begin{bmatrix}\ \lambda \!&\! 0\ \\ \ -\lambda \!&\! \lambda\ 
\end{bmatrix}
\begin{bmatrix}\ x\ \\ \ y\ \end{bmatrix}
+\begin{bmatrix}\ p(x,y)\ \\ \ q(x,y)\ \end{bmatrix}.
\end{align*}
Also, we transform it to a form of polar coordinates:
\begin{align}
\begin{aligned}
r'(t)&=\lambda r-\lambda r\sin\theta\cos\theta+r\varphi(r,\theta),\quad 
\theta'(t)=-\lambda\cos^2\theta+\psi(r,\theta).
\end{aligned}
\label{eqn:10}
\end{align}
Since both $\varphi,\psi$ are $\alpha$-H\"older continuous, 
there are reals $c>0,m>0$ such that 
\begin{align}
\big|\,\varphi(r,\theta)\,\big|<mr^{\alpha},\quad
\big|\,\psi(r,\theta)\,\big|<mr^{\alpha}\quad
\textrm{on }[0,c]\times\mathbb{R}.
\label{eqn:11}
\end{align}
We take a real $d$ so that $md^{\alpha}<\tfrac{1}{4}|\lambda|$ 
and $d\leqq c$. 
If we take $r(t)\in[0,d]$, then we evaluate the first of 
Eq.(\ref{eqn:10}) as follows:
\begin{align*}
r'(t)\leqq\big(\!-\!|\lambda|+\tfrac{1}{2}|\lambda|
+\tfrac{1}{4}|\lambda|\,\big)r=-kr,\quad 
k:=\tfrac{1}{4}|\lambda|. 
\end{align*}
If we take $r_0=r(0)$ so small that $r_0\in[0,d]$, 
then we have that
\begin{align}
r(t)\leqq r_0e^{-kt}.
\label{eqn:12}
\end{align}
Therefore, $r(t)\to 0$ as $t\to\infty$. 
Hence $\bm{p}$ is attracting. \par
\medskip
By applying Eqs.(\ref{eqn:11}), (\ref{eqn:12}) to the second 
of Eq.(\ref{eqn:10}), we have that
\begin{align}
a\cos^2\theta-be^{-lt}\leqq\theta'(t)\leqq 
a\cos^2\theta+be^{-lt},\quad& 
\label{eqn:13} \\
a:=|\lambda|,\ b:=mr_0^{\,\alpha},\ l:=\alpha k.&
\notag
\end{align}
Take a real $t_0\in[0,\infty)$ so large that 
$e^{lt_0}>16ab/l^2$. \par
\medskip
Firstly, consider the case when $h(t)<0$ for every $t\in[t_0,\infty)$, 
where
\begin{align*}
h(t)=a\cos^2\theta(t)-be^{-lt}.
\end{align*}
In this case, we have that $\cos\theta(t)\to 0$ as $t\to\infty$. 
Therefore, there is an integer $n$ such that 
$\theta(t)\to (n+\tfrac{1}{2})\pi$ as $t\to\infty$. 
Hence $\bm{p}$ is a node. \par
\medskip
Secondly, consider 
the case when $h(t_1)\geqq 0$ for some $t_1\in[t_0,\infty)$. 
In this case, we shall prove that $h(t)\geqq 0$ for every 
$t\in[t_1,\infty)$. 
To prove it, assume that there is a real $t_2\in(t_1,\infty)$ 
satisfying $h(t_2)<0$. 
Take the largest $s\in[t_1,t_2)$ satisfying $h(s)=0$. 
This implies that $h(t)<0$ on $(s,t_2]$. 
By the definition of derivative, we have that
\begin{align}
h'(s)=\lim_{t\to s+0}\dfrac{h(t)-h(s)}{t-s}\leqq 0. 
\label{eqn:14}
\end{align}
By differentiating $h(t)$, and by putting $t=s$, we have that
\begin{align}
h'(s)=-2a\cos\theta(s)\sin\theta(s)\,\theta'(s)+lbe^{-ls}.
\label{eqn:15}
\end{align}
By putting $h(s)=0$ to Eq.(\ref{eqn:13}), we have that 
\begin{align}
0\leqq\theta'(s)\leqq a\cos^2\theta(s)+be^{-ls}.
\label{eqn:16}
\end{align}
By using Eqs.(\ref{eqn:14}), (\ref{eqn:15}), we have that
\begin{align}
\cos\theta(s)\sin\theta(s)\,\theta'(s)
=\dfrac{1}{2a}\big\{lbe^{-ls}-h'(s)\big\}>0. 
\label{eqn:17}
\end{align}
By applying Eq.(\ref{eqn:16}) to Eq.(\ref{eqn:17}), we have that
\begin{align}
\cos\theta(s)\sin\theta(s)>0,\quad\theta'(s)>0.
\label{eqn:18}
\end{align}
By applying Eqs.(\ref{eqn:16}), (\ref{eqn:18}) to Eq.(\ref{eqn:15}), 
we have that
\begin{align}
h'(s)\geqq-2a\cos\theta(s)\sin\theta(s)
\big(a\cos^2\theta(s)+be^{-ls}\big)+lbe^{-ls}.
\label{eqn:19}
\end{align}
Since $a\cos^2\theta(s)=be^{-ls}$, we have that
\begin{align}
&\hspace{-2em}
\big(lbe^{-ls}\big)^2-\big\{2a\cos\theta(s)\sin\theta(s)
\big(a\cos^2\theta(s)+be^{-ls}\big)\big\}^2 \notag\\
&=\big(lbe^{-ls}\big)^2-4a^2\dfrac{be^{-ls}}{a}
\Big(1-\dfrac{be^{-ls}}{a}\Big)\big(be^{-ls}+be^{-ls}\big)^2 \notag\\
&=l^2b^2e^{-3ls}\big(e^{ls}-16ab/l^2\big)+16b^4e^{-4ls}.
\label{eqn:20}
\end{align}
Since $e^{ls}\geqq e^{lt_0}>16ab/l^2$, 
we have that Eq.(\ref{eqn:20})\ $>0$. 
By applying it to Eq.(\ref{eqn:19}), we have that $h'(s)>0$. 
It contradicts to Eq.(\ref{eqn:13}). 
Therefore, we have that  $h(t)\geqq 0$ for every $t\in[t_1,\infty)$. \par
\medskip
By using Eq.(\ref{eqn:13}), we have that $\theta(t)$ is monotone 
increasing on $[t_1,\infty)$. 
Moreover, $\theta(t)$ can not cross $(n+\tfrac{1}{2})\pi$ 
for any integer $n$. 
This is because if $\theta(t_2)=(n+\tfrac{1}{2})\pi$ for some real 
$t_2\in[t_1\infty)$, then $h(t_2)=-be^{-kt_2}<0$, a contradiction. 
Therefore, $\theta(t)$ is bounded. 
Since $\theta(t)$ is monotone increasing and bounded, it converges 
to a limit. Hence $\bm{p}$ is a node. 
\end{proof}
\end{document}